\documentclass[12pt]{article}
\usepackage[T2A]{fontenc}
\usepackage[cp866]{inputenc}
\usepackage[russian]{babel}

\usepackage{amssymb}
\usepackage{amsmath}
\usepackage{euscript}

\begin{document}

\centerline{\large\bf Finite-dimensional reduction of systems}

\centerline{\large\bf of nonlinear diffusion equations }

\bigskip
\centerline{\large\bf A.V. Romanov}
\bigskip

\textbf{Abstract.} We present a class of one-dimensional systems
of nonlinear parabolic equations for which long-time phase
dynamics can be described by an ODE with a Lipschitz vector field
in $\mathbb{R}^{n}$. In the considered case of the Dirichlet
boundary value problem sufficient conditions for a
finite-dimensional reduction turn out to be much wider than the
known conditions of this kind for a periodic situation.

\bigskip
{\it Keywords}: nonlinear parabolic equations, finite-dimensional dynamics on attractor;
inertial manifold.

\textit{2020 Mathematics Subject Classification$:$ Primary 35B41,
35K57; Secondary 35K42, 35K90, 35K91}

\bigskip
\noindent{\bf 1. Introduction}
\medskip

One of the main problems in the study of evolution equations
is related to describing the final (at large time) behavior
of their solutions.
We consider systems of diffusion equations with Dirichlet boundary condition
$$
\partial_{t}u=D\partial_{xx}u+f(x,u)\partial_{x}u+g(x,u),\quad \quad
u(0)=u(1)=0
\eqno (1.1)
$$
on the closed interval $J=[0,1]$. Here $u=(u_{1},\dots,u_{m})$,
$f$ and $g$ are sufficiently regular matrix and vector functions,
respectively. We assume that the matrix $D$ of numerical
coefficients is similar to a diagonal matrix with positive
eigenvalues. In the case of
$D=\mathrm{diag}\,\{d_{1},\dots,d_{m}\}$, $d_{j}>0$, we deal with
reaction-diffusion-convection equations. Under appropriate
conditions on $f$ and $g$, system (1.1) induces a smooth
dissipative semiflow $\{\, \Phi_{t} \}_{t\geq 0}$ in the phase
space $X^{\alpha}\subset C^{1}(J,\mathbb{R}^{m})$ with an
appropriate $\alpha>0$, where $\{X^{\alpha}\}_{\alpha\geq 0}$ is
the Hilbert semiscale [3] generated by the linear sectorial
operator $u\rightarrow -Du_{xx}$ in $X=L^{2}(J,\mathbb{R}^{m})$.
In this situation, there exists a global attractor [2,7,12] (in
what follows, simply an attractor), i.e., a connected compact
invariant set $\mathcal{A}\subset X^{\alpha}$ of a finite
Hausdorff dimension uniformly attracting bounded subsets of
$X^{\alpha}$ as $t\to +\infty $.

Our goal is to find conditions under which the dynamics on the
attractor (final dynamics) of parabolic system (1.1) is
finite-dimensional in the sense of [8]. This means that, for some
ODE $\partial_{t}\xi=h(\xi)$ in ${\mathbb R}^{N}$ with Lipschitz
vector field $h$, the resolving flow $\{ \, \Theta_{t} \}$ and an
invariant compact set ${\EuScript K}\subset {\mathbb R}^{N}$, the
phase semiflows $\{\, \Phi_{t} \}_{t\geq 0}$ on ${\EuScript A}\;$
and $\{ \, \Theta_{t} \}_{t\geq 0}$ on ${\EuScript K}$ are
Lipschitz-conjugate. In this connection, we can speak [14] about
the finite-dimensional reduction of evolution problem (1.1).

The main result in this paper (Theorem~4.3) ensures
that the final phase dynamics of system (1.1)
is finite-dimensional  under the \textit{consistency condition}
$$
Df(x,u)=f(x,u)D,\quad \quad (x,u)\in {J\times\rm co}\,\mathcal{A}, \eqno(1.2)
$$
where $\mathrm{co}\,\mathcal{A}$ is the convex hull of $\mathcal{A}$.

It is well known [5] that, in the case of scalar diffusion ($D=dE$
with unit matrix $E$) and $f=f(u)$, $g=g(u)$, there exists an
inertial manifold (IM), i.e., a finite-dimensional invariant
$C^{1}$-surface in the phase space containing an attractor and
exponentially attracting (with an asymptotic phase) all
trajectories of the system as $t\rightarrow +\infty$. The presence
of IM implies that the final dynamics is finite-dimensional, and
an extensive literature deals with the existence of such manifolds
(see, [7,10,12,14]). An original approach to these problems is
presented in recent works of M.~Anikushin (see [1] and the
references therein).

In the periodic case ($J$ is a circle of length~1), conditions
ensuring that the final dynamics of systems (1.1) with
$D=\mathrm{diag}$ is finite-dimensional were obtained by the
author in [11; p.13409]. Note that, in the class of periodic
systems (1.1) with scalar diffusion, the first example of
semilinear parabolic equation of mathematical physics that does
not demonstrate such a dynamics was constructed in [6;
Theorem~1.2].

\bigskip
\noindent{\bf 2. Preliminaries}
\medskip

In what follows, if necessary, we will use the technique developed
in~[11]. All preliminary constructions in Sections~2 and~3 are
carried out for the case $D=\mathrm{diag}$. Let us write system
(1.1) in the form of a semilinear parabolic equation (SPE)
$$
\partial _{t} u=-Au+F(u)
\eqno (2.1)
$$
in the \textit{real} Hilbert space $X=L^{2}(J,\mathbb{R}^{m})$ equipped with the norm $\|\cdot\|$.
Here we have $A:u\rightarrow -Du_{xx}$ with Dirichlet boundary condition
and the nonlinearity $F:u\rightarrow f(x,u)\partial_{x}u+g(u)$.
For the linear positive definite operator $A$, we put
$X^{\alpha} =\mathcal{D}(A^{\alpha})$ with $\alpha \geq 0$ and $X_{0}=X$.
Then $\left\|u\right\|_{\alpha}=\left\|A^{\alpha}u\right\|$.
Note that the function $F$ is of class $W^{2}(X^{\alpha},X)$ if
$$
F\in C^{2}(X^{\alpha},X) \bigcap \mathrm{Lip}(X^{\alpha},X)
\quad\text{and}\quad
\|F(u)\|\leq M,\quad u\in X^{\alpha},
\eqno (2.2)
$$
for some $\alpha \in [0,1)$. In this case, SPE (2.1) generates [3]
a smooth compact resolving semiflow $\{ \Phi_{t} \}_{t\ge 0}$ in
the phase space $X^{\alpha }$. Assumption (2.2) implies [10;
Lemma~1.1] the $X^{\alpha}$-dissipativity of (2.1):
$$
\mathop{\limsup}\limits_{t\rightarrow+\infty}\,\left\|\Phi_{t}u\right\|_{\alpha} \le r
$$
for some $r>0$ uniformly in $u\in$ balls in $X^{\alpha}$. Under
such conditions, there exists [2,7,12] a compact attractor
$\mathcal{A}\subset X^{\alpha}$ consisting of all bounded complete
trajectories $\{u(t)\}_{t\in \mathbb{R}} \subset X^{\alpha }$. In
fact, $\mathcal{A}\subset X^{1}$ due to the \textit{smoothing
action} of the parabolic equation [3]. Simple argument [11;
p.13410] shows that, in all constructions concerning SPE (2.1),
one can replace the nonlinearity exponent $\alpha$ by any value
$\alpha_{1}\in(\alpha,1)$, and if condition (2.2) is satisfied in
the pair of spaces $(X^{\theta},X^{\theta+\alpha})$ with $\theta
>0$ instead of $(X,X^{\alpha})$, then all the listed properties of
the dynamics are preserved for the phase space
$X^{\theta+\alpha}$. In what follows, functions $Y_{1}\rightarrow
Y_{2}$ of class (2.2) will arise for some Banach spaces
$Y_{1},Y_{2}$.

As in [11], we will use sufficient conditions for the final
dynamics to be finite-dimensional~[9]. Assume that $G(u)=F(u)-Au$
is the vector field (2.1),
$\mathcal{N}=\mathcal{A}\times\mathcal{A}$, and $Y$ is a Banach
space.
\medskip

\textbf{Definition~2.1 ([9]).} A continuous field
$\Pi:\mathcal{N}\rightarrow Y$ is said to be regular if, for any
$u,v\in \mathcal{A}$, the function
$\Pi(\Phi_{t}u,\Phi_{t}v):[0,+\infty)\rightarrow Y$ is of class
$C^{1}$ with the derivative $\partial_{t}\Pi(u,v)$ uniformly
bounded in $(u,v)\in \mathcal{N}$ at zero.
\medskip

The smoothness of the semiflow $\{\Phi_{t}\}$ and the invariance
of the compact set $\mathcal{A}\subset X^{\alpha}$ ensure the
regularity of the identity embedding $\mathcal{N}\rightarrow
X^{\alpha}\times X^{\alpha}$, and hence, the regularity of each
field $\Pi:\mathcal{N}\rightarrow Y$ that can be continued to a
$C^{1}$-mapping into the $X^{\alpha}\times
X^{\alpha}$-neighborhood of the set $\mathcal{N}$. In this
situation, we have $\partial_{t}\Pi(u,v)=\Pi'(u,v)(G(u),G(v))$,
where $(\,\cdot\,)'$ is the Frechet differentiation. Under
condition (2.2) on the nonlinearity $F$, the function
$u\rightarrow G(u)$ on $\mathcal{A}$ is continuous and even
H\"older [8] in the $X^{\alpha}$-metric. The regular fields
$\mathcal{N}\rightarrow Y$ form a linear structure and even a
multiplicative one if $Y$ is a Banach algebra. In the latter case,
if all elements $\Pi(u,v)\in Y$ are invertible, then the field
$\Pi^{-1}$ is also regular.

We will start from the decomposition
$$
G(u)-G(v)=(T_{0}(u,v)-T(u,v))(u-v), \quad (u,v) \in \mathcal{\mathcal{N}},
\eqno (2.3)
$$
where $T_{0}\in \mathcal{L}(X^{\alpha})$ and $T\in\mathcal{L}(X^{1},X)$
are unbounded linear operators in $X$
similar to positive definite ones.
We let
$$
\Sigma_{T} =\bigcup_{u,v\in \mathcal{A}} \mathrm{spec}\,T(u,v)
$$
denote the total spectrum of the operators $T$.

We will need a particular case [9; Theorem~2.8] in the situation
$\Sigma_{T}\subset \mathbb{R}^{+}$.
\medskip

\textbf{Theorem 2.2.}
\textit{Assume that $F\in W^{2}(X^{\alpha},X)$ and
$$
T(u,v)= S^{-1}(u,v)H(u,v)S(u,v)
\eqno (2.4)
$$
on $\mathcal{N}$,
where the unbounded self-adjoint linear operators $H(u,v)$
are positive definite in $X$,
the fields
$S,\,S^{-1}:\mathcal{N}\rightarrow \mathcal{L}(X)$ and
$T_{0}:\mathcal{N}\rightarrow \mathcal{L}(X^{\alpha},X)$
are regular, and the field
$T_{0}:\mathcal{N}\rightarrow\mathcal{L}(X^{\alpha})$
is bounded. If in addition, the set $\mathbb{R}^{+}\backslash \Sigma_{T}$
contains intervals $(a_{k}-\xi_{k},a_{k}+\xi_{k})$ with $a_{k}>\xi_{k}>0$
such that
$$
\xi_{k}\rightarrow \infty,\quad\quad a_{k}^{\alpha/2}=o(\xi_{k})
\eqno (2.5)
$$
as $k \to +\infty$, then the final $X^{\alpha}$-dynamics of SPE $(2.1)$
is finite-dimensional.}
\medskip

We further assume that the matrix function $f=f(x,u)$ and the vector functions
$g=g(x,u)$ in (1.1) satisfy the regularity conditions:

\textbf{(H)} $\quad f,g$ of class $C^{\infty}$ on $J\times {\mathbb R}^{m}$
are finite in $u$ and $f(x,0)=g(x,0)=0$ for $x=0,1$.

\medskip
We let $\mathcal{H}^{s}=\mathcal{H}^{s}(J)$ denote generalized
Sobolev $L^{2}$-spaces (spaces of Bessel potentials [3,13]) of
scalar functions on $J$ with arbitrary $s\geq 0$. If $s>1/2$, then
$\mathcal{H}^{s}\subset C(J)$ and $\mathcal{H}^{s}$ is a Banach
algebra [13; Sec.~2.8.3]. The differentiation operator
$\partial_{x}$ belongs to
$\mathcal{L}(\mathcal{H}^{s+1},\mathcal{H}^{s})$. In fact, the
$X^{s}$ are closed subspaces (with equivalent norm) in the spaces
$\mathcal{H}^{2s}(J,\mathbb{R}^{m})$ of vector functions, and
$X^{s}=\mathcal{H}^{2s}(J,\mathbb{R}^{m})$ for $s\leq 1/4$. For
$s>1/4$, the space $X^{s}$ consists of elements $u\in
\mathcal{H}^{2s}(J,\mathbb{R}^{m})$ such that $u(0)=u(1)=0$.

We now fix an arbitrary $\alpha \in (3/4,1)$. Then we have
$\mathcal{H}^{2\alpha}\hookrightarrow C^{1}(J)$ and
$X^{\alpha}\hookrightarrow C^{1}(J)$, where the symbol
$\hookrightarrow$ denotes a linear continuous embedding of
function spaces. We will use several required embedding theorems
[3,13]. For an arbitrary $C^{\infty}$-function $z:J\times
\mathbb{R}^{m}\rightarrow \mathbb{R}$, the mapping
$\psi:u\rightarrow z(x,u)$ is a function of class $W^{2}$ (see
(2.2)) from $C^{s}(J,\mathbb{R}^{m})$ to $C^{s}(J)$ for all
$s\in\mathbb{N}$. This implies that $\psi\in
W^{2}(\mathcal{H}^{2\alpha}(J,\mathbb{R}^{m}),C^{1}(J))$. Using
the embeddings $\mathcal{H}^{s+1}\hookrightarrow
C^{s}(J)\hookrightarrow \mathcal{H}^{s}$, we can conclude that
$\psi\in
W^{2}(\mathcal{H}^{s}(J,\mathbb{R}^{m}),\mathcal{H}^{s}(J))$. We
thus obtain $\,F\in W^{2}(X^{1},X^{1/2})$ for the nonlinear part
$F:u\rightarrow f(x,u)\partial_{x}u+g(u)$ of system (1.1).
Moreover, $X^{\alpha}\hookrightarrow
C^{1}(J,\mathbb{R}^{m})\hookrightarrow
C(J,\mathbb{R}^{m})\hookrightarrow X$, and hence $F\in
W^{2}(X^{\alpha},X)$. We also note that $X^{3/2}\hookrightarrow
C^{2}(J,\mathbb{R}^{m})$.

We choose $X^{\alpha}$ as the phase space of system (1.1).
Then the phase dynamics of (1.1) in $X^{\alpha}$
is dissipative, and there exists a global attractor
$\mathcal{A} \subset X^{\alpha}$.
Since $\,F\in W^{2}(X^{1},X^{1/2})$, system (1.1) also
generates a smooth dissipative phase semiflow in the space $X^{1}$
and the attractor $\mathcal{A}$ is compact in $X^{3/2}$.
As above, we denote $\mathcal{N}=\mathcal{A}\times \mathcal{A}$.
\medskip

\textbf{Remark 2.3}. The phase dynamics of system (1.1) has the
following property: if $Y$ is a Banach space, then each vector
field $\Pi:\mathcal{N}\rightarrow Y$ continuous in the
$(X^{\alpha}\times X^{\alpha})$-metric and extendable to a
$C^{1}$-mapping $X^{1}\times X^{1}\rightarrow Y$ is regular in the
sense of Definition~2.1.

Indeed, the smoothness of the semiflow in $X^{1}$ means that the mapping
$(t,u)\rightarrow \Phi_{t}u:\,(0,+\infty)\times X^{1}\rightarrow X^{1}$ is smooth.
This ensures that the identity embedding $\mathcal{N}\rightarrow X^{1}\times X^{1}$
is regular, and hence the field $\Pi$ on $\mathcal{N}$ is also regular.

\bigskip
\noindent{\bf 3. Decomposition of a vector field on an attractor}
\medskip

We want to apply Theorem~2.2 to SPE (1.1) with $D=\mathrm{diag}$
and the phase space $X^{\alpha}$, $\alpha\in (3/4,1)$. We let
$\mathbb{M}^{m}$ denote the algebra of numerical $m\times m$
matrices with Euclidean norm, and let $Y(J,\mathbb{M}^{m})$ denote
linear spaces of such matrices with elements from some Banach
space $Y$ of scalar functions on $J=[0,1]$. Similarly [11;
pp.13412--13413 ], we assume that
$$
B_{0}(x;u,v)=\int_{0}^{1}(f_{u}(x,w(x))w_{x}(x)+g_{u}(x,w(x))d\tau,
\eqno (3.1.1)
$$
$$
B(x;u,v)=\int_{0}^{1}f(x,w(x))d\tau
\eqno (3.1.2)
$$
for $u,v\in X^{\alpha},\; w(x)=\tau u(x)+(1-\tau)v(x)$, $x\in J$.
The elements of the matrices $B_{0},B$ are continuous functions,
and for $u,v\in \mathcal{A}$, they are functions of class $C^{2}$ on $J$.
Using the $C^{1}$-smoothness of the mappings
$(u,v)\rightarrow f_{u}(x,w)w_{x}+g_{u}(x,w),\; (u,v)\rightarrow
f(x,w),\;\, X^{\alpha}\times X^{\alpha}\rightarrow C(J,\mathbb{M}^{m})$
for a fixed $\tau\in [0,1]$ and differentiating the integrands in the expressions for $B_{0}$ and $B$
with respect to the parameters $(u,v)$, we see that the mappings
$$
(u,v)\rightarrow B_{0}(\cdot\,;u,v),\quad (u,v)\rightarrow B(\cdot\,;u,v) \eqno(3.2)
$$
are of class $C^{1}(X^{\alpha}\times X^{\alpha},C(J,\mathbb{M}^{m}))$.
We use the integral mean value theorem for nonlinear operators to write
the decomposition of the vector field of (1.1) on the  attractor
$\mathcal{A}\subset X^{\alpha}$ as
$$
G(u)-G(v)=-Ah+(\int_{0}^{1}F'(\tau u+(1-\tau)v)d\tau)h
$$
$$
=Dh_{xx}+B_{0}(x;u,v)h+B(x;u,v)h_{x}, \quad u,v\in \mathcal{A},
$$
where $h=u-v$, $\tau u+(1-\tau)v\in {\rm co}\,\mathcal{A}$, and
$(\,\cdot\,)'$ is the Frechet differentiation. To eliminate the
dependence on $h_{x}$, we apply (following [4]) the transformation
$h=U\eta$, where the $m\times m$ matrix function $U(x)=U(x;u,v)$,
$x\in [0,1]$, is the solution of the linear Cauchy problem
$$
U_{x}=-\frac{1}{2} D^{-1} B(x)U,\quad \; U(0)=E.
\eqno (3.3)
$$
As a result, we obtain relation (2.3)  with linear operators
$$
T_{0} (u,v)h=(B_{0}(x)-\frac{1}{2} B_{x} (x)-\frac{1}{4} B(x)D^{-1} B(x))h,
\eqno (3.4.1)
$$
$$
T(u,v)h=-DU\partial _{xx} U^{-1} h.
\eqno(3.4.2)
$$

We note that the change of variable $h=U\eta$
does not change the Dirichlet boundary conditions
for the linear part of (1.1).
In the expressions for the matrices $B_{0}$, $B$, $U$, $U^{-1}$,
we often omit the dependence on $u$, $v$, and sometimes on $x$.
\medskip

\textbf{Lemma 3.1.}
\textit{The field of operators $T_{0}$ on $\mathcal{N}$
is regular with values in $\mathcal{L}(X^{\alpha},X)$
and bounded with values in $\mathcal{L}(X^{\alpha})$.}
\medskip

\textbf{Proof.}
We assume that $T_{0}h=Q(x;u,v)h$ in (3.4.1) with
$h\in\mathcal{A}-\mathcal{A} \subset X^{\alpha}$.
The convex hull of the attractor $\mathcal{A}$ is bounded in the $X^{3/2}$-norm
equivalent to the $\mathcal{H}^{3}(J,\mathbb{R}^{m})$-norm,
and hence the matrix functions $B$, $BD^{-1}B$ and $B_{0}$ are bounded
uniformly in $(u,v)\in \mathcal{N}$ in $\mathcal{H}^{3}(J,\mathbb{M}^{m})$
and $\mathcal{H}^{2}(J,\mathbb{M}^{m})$, respectively.
Thus, the matrix functions $B_{x}$ and $Q$ are bounded on $\mathcal{N}$
in the norm of $\mathcal{H}^{2}(J,\mathbb{M}^{m})$
and $T_{0}$ is the operator of multiplication of vector functions
from $X^{\alpha}\subset \mathcal{H}^{2\alpha}(J,\mathbb{R}^{m})$
by the matrix $Q\in \mathcal{H}^{2\alpha}(J,\mathbb{M}^{m})$
with $2\alpha \in (3/2,2)$.
Since $\mathcal{H}^{2\alpha}(J)$ is a Banach algebra,
we obtain $T_{0}(u,v)\in \mathcal{L}(X^{\alpha})$
and $\|T_{0}(u,v)\|_{\alpha}\leq \mathrm{const}$ on $\mathcal{N}$.

With regard to Remark~2.3 and the above-noted smoothness of
mappings (3.2), the regularity of the field of the operators
$T_{0}:\mathcal{N}\rightarrow\mathcal{L}(X^{\alpha},X)$ can be
proved as in the case of periodic boundary conditions in [11;
Lemma~3.3]. $\;\Box$
\medskip

The matrix function $U(x)$ in the Cauchy problem (3.3) can be treated as
a bounded linear operator in $X$.
\medskip

\textbf{Lemma 3.2.}
\textit{The fields of the operators
$U,U^{-1}:\mathcal{N}\rightarrow \mathcal{L}(X)$ are regular.}
\medskip

For the field of $U$, this can be proved as a similar assertion in
the periodic case [11; Lemma~3.4]. At the same time, the
regularity of $U$ implies the regularity of the field of the
inverse operators $U^{-1}$.

\medskip
Now we assume that $d_{-}=\mathop{\min }\limits_{1\leq j\leq m}\,d_{j}$
and $d_{+}=\mathop{\max }\limits_{1\leq j\leq m}\,d_{j}$ for $D=\mathrm{diag}\{d_{j}\}$.
Assume also that $\{\lambda_{n}:\lambda_{1}<\lambda_{2}<\dots\}$
are eigenvalues of the linear operator $A=-D\partial_{xx}$.
Since
$$
\mathrm{spec}\,A=\{d_{j}\pi ^{2} \nu ^{2},\;\; \nu \in {\mathbb
N},\; \; j\in \overline{1,m}\,\}, \eqno(3.5)
$$
we have $\lambda_{n}\leq\pi^{2}d_{+} \,n^{2}$. Using the counting
function for $\mathrm{spec}\,A$, we obtain
$$
n\leq\sum_{j=1}^{m}\dfrac{\sqrt{\lambda_{n}}}{\pi \sqrt{d_{j}}}\leq \frac {m} {\pi
\sqrt{d_{-}}}\,\sqrt{\lambda_{n}}\,,
$$
and hence
$$
\frac {\pi^{2}d_{-}} {m^{2}} \,n^{2}\leq \lambda_{n}\leq\pi^{2}d_{+}
\,n^{2},\quad n\in\mathbb{N}.
\eqno (3.6)
$$

\medskip
\textbf{Lemma 3.3.}
\textit{The estimate
$\mathop{\limsup}\limits_{n\rightarrow \infty}\,
n^{-1}(\lambda_{n+1}-\lambda_{n})>0$ holds.}

\medskip
{\bf Proof.}
If, on the contrary,
$\lambda_{n+1}-\lambda_{n}=\beta_{n}n$ with
$\beta_{n}\stackrel{n\rightarrow\infty}{\longrightarrow}0$,
then
$$
n^{-2}\lambda_{n}=n^{-2}(\lambda_{1}+\mathop{\sum}\limits_{k=1}\limits^{n-1}(\lambda_{k+1}-\lambda_{k})=
n^{-2}(\lambda_{1}+\mathop{\sum}\limits_{k=1}\limits^{n-1}\beta_{k}k)
$$
$$
\leq
n^{-2}(\lambda_{1}+\mathop{\sum}\limits_{k=1}\limits^{n-1}\beta_{k}n)\leq
n^{-2}\lambda_{1}+n^{-1}
\mathop{\sum}\limits_{k=1}\limits^{n}\beta_{k}.
$$
However, this implies the relation $\lambda_{n}=o(n^{2})$
which contradicts the left inequality in (3.6). $\;\Box$

\bigskip
\noindent{\bf 4. Main results}
\medskip

By the assumptions of Theorem~2.1, it is necessary to prove
the ``uniform'' similarity of the operators $T(u,v)$ in (3.4.2)
to positive definite operators of the form (2.4),
as well as the required sparsity (2.5) of their total spectrum $\Sigma_{T}$.
We assume that the regularity conditions \textbf{(H)}
are satisfied for the functions $f$ and $g$ in (1.1).

\medskip
\textbf{Theorem 4.1.}
\textit{If the matrix $D=\mathrm{diag}\{d_{j}\}$ with $d_{j}>0$
and condition $(1.2)$ is satisfied, then the phase dynamics on the attractor
is finite-dimensional}.
\medskip

\textbf{Proof.}
The operator $A=-D\partial_{xx}$ with Dirichlet condition is self-adjoint
and positive definite in $X$.
Assumption (1.2) (for any $x\in J$ and $u,v\in \mathcal{A}$)
implies the relation $DB(x)=B(x)D$ for the matrices $B(x)=B(x;u,v)$ in (3.1.2).
Thus, the matrices $B(x)$ and $D^{-1}B(x)$ inherit the block
(with respect to equal $d_{j}$) structure of the diffusion matrix
$D=\mathrm{diag}\,\{d_{1},\dots,d_{m}\}$.
Therefore, the same also holds for the solutions $U(x)$ of the Cauchy problem (3.3),
and hence, $DU(x)=U(x)D,\;x\in J$, and
$$
T(u,v)=U(u,v)(-D\partial_{xx})U^{-1}(u,v)
$$
in (3.4.2).
Thus, representation (2.4) with $S(u,v)=U^{-1} (u,v)$ and $H(u,v)\equiv A$
holds for $T(u,v)$.
The total spectrum $\Sigma_{T}$ coincides with $\mathrm{spec}(A)$ in (3.5).
By Lemma~3.3, there exists $\varepsilon>0$ and an increasing sequence
of indices $n(k)$ such that $\lambda_{n(k)+1}-\lambda_{n(k)}>\varepsilon n(k)$
for $k\geq k_{0}$.
We put $a_{k}=(\lambda_{n(k)+1}+\lambda_{n(k)})/2$,
$\xi_{k}=(\lambda_{n(k)+1}-\lambda_{n(k)})/3$ and $M=\pi^{2}d_{+}$.
From the right inequality in (3.6) we obtain
$$
a_{k}\leq M(n^{2}(k)+n(k)+\frac {1} {2})\leq 3Mn^{2}(k)\leq \frac {3M}
{\varepsilon^{2}}(\lambda_{n(k)+1}-\lambda_{n(k)})^{2}\leq \frac
{27M} {\varepsilon^{2}}\xi_{k}^{2}
$$
for $k\geq k_{0}$, i.e., $a_{k}=O(\xi_{k}^{2})$ as $k\rightarrow \infty$.
Since $a_{k}^{\alpha/2}=o(\xi_{k})$ for $\alpha\in (3/4,1)$ and $k\rightarrow\infty$,
the sought assertion follows from Lemmas~3.1 and~3.2 and Theorem~2.2. $\Box$
\medskip

\textbf{Remark 4.2.}
Parabolic systems (1.1) with $D=\mathrm{diag}$ demonstrate a finite-dimensional dynamics
on the attractor for any admissible nonlinearities $f$ and $g$
in the case of scalar diffusion
and under the condition $f=\mathrm{diag}$ in the case of $m$
distinct diffusion coefficients $d_{j}$.
In the case of $s$ distinct diffusion coefficients with $1<s<m$,
the dynamics on the attractor is finite-dimensional under the condition that
the matrix function $f$ inherits the block (with respect to the same $d_{j}$)
structure of the matrix $D=\mathrm{diag}\{d_{j}\}$.

\medskip
Now we formulate the main result.
We assume that the matrix $D$ in system (1.1) has the form $D=C\overline{D}C^{-1}$,
where the matrix $C$ is nondegenerate and
$\overline{D}=\mathrm{diag}\,\{d_{1},\dots,d_{m}\}$ with $d_{j}>0$.
The linear operator $-D\partial_{xx}=-C(\overline{D}\partial_{xx})C^{-1}$
is sectorial in $X=L^{2}(J,\mathbb{R}^{m})$.
The change of variable $u=Cv$ reduces (1.1) to the system of equations
$$
\partial_{t}v=\overline{D}\partial_{xx}v+\overline{f}(x,v)\partial_{x}v+\overline{g}(x,v),\quad
\quad v(0)=v(1)=0,
$$
$$
\overline{f}(x,v)=C^{-1}f(x,Cv)C, \quad \; \overline{g}(x,v)=C^{-1}{g}(x,Cv).
\eqno(4.2)
$$
The matrix function $\overline{f}$ and the vector function $\overline{g}$ inherit
the regularity properties \textbf{(H)} of the original functions $f$ and $g$.
The system of equations (4.2) is dissipative in $X^{\alpha}$,
and hence, the same is also true for system (1.1).
The attractors $\mathcal{A}$ of system (1.1) and $\overline{\mathcal{A}}$ of system (4.2)
are related by the formula $\mathcal{A}=C\overline{\mathcal{A}}$.
By the definition of the finite-dimensionality of the final phase dynamics
(Section~1), systems (4.2) and (1.1) simultaneously demonstrate this property.
\medskip

\textbf{Theorem 4.3 (main theorem).}
\textit{If the matrix $D$ is similar to $\mathrm{diag}\{d_{j}\}$ for $d_{j}>0$
and consistency condition $(1.2)$ is satisfied, then the final dynamics
of system $(1.1)$ is finite-dimensional}.

\medskip
\textbf{Proof.}
Since $Df(x,u)=f(x,u)D$ on $J\times\rm {co}\,\mathcal{A}$, we have
$$
\overline{D}\,\overline{f}(x,v)=C^{-1}DC\cdot C^{-1}f(x,Cv)C=C^{-1}Df(x,u)C
$$
$$
=C^{-1}f(x,u)DC=C^{-1}f(x,Cv)C\overline{D}=\overline{f}(x,v)\overline {D}
$$
on $J\times \rm {co}\,\overline{\mathcal{A}}$.
Here $u\in \rm {co}\,\mathcal{A}$ and $v\in \overline{\mathcal{A}}$.
So we see that condition (1.2) is satisfied for the matrix function $\overline{f}$
and, by Theorem~4.1, the dynamics of system (4.2) on the attractor
$\overline{\mathcal{A}}\subset X^{\alpha}$ is finite-dimensional.
This also implies that the dynamics of system (1.1) is finite-dimensional
on the attractor $\mathcal{A}\subset X^{\alpha}$. $\;\Box$

\medskip
\textbf{Remark 4.4.}
Under consistency condition (1.2), the final dynamics of system (1.1) is finite-dimensional
if all eigenvalues of the matrix $D$ are distinct and positive.
In particular, condition (1.2) is satisfied for $f=D_{1}\varphi$,
where the numerical matrix $D_{1}$ commutes with $D$
and $\varphi=\varphi(x,u)$ is a smooth scalar function finite in $u$.

\bigskip
\noindent \textbf{References}
\medskip

\noindent [1] M. Anikushin, ``Frequency theorem for parabolic
equations and its relation to

inertial manifolds theory``, \textit{J. Math. Anal. and Appl.},
\textbf{505}:1 125454 (2022).

\noindent [2] L. V. Babin, M. I. Vishik, \textit{Attractors of
evolution equations}, North-Holland, 1992.

\noindent [3] D. Henry, \textit{Geometric theory of semilinear
parabolic equations}, Lect. Notes in Math.,

vol. \textbf{840}, Springer, 1981.

\noindent [4] D. A. Kamaev, ``Families of stable manifolds of
invariant sets of systems of parabolic

equations'', \textit{Russ. Math. Surv.}, \textbf{47}:5 (1992),
185-186.

\noindent [5] A. Kostianko and S. Zelik, ``Inertial manifolds for
1D reaction-diffusion-advection

systems. Part I: Dirichlet and Neumann boundary conditions'',
\textit{Comm. Pure Appl.}

\textit{Anal.,} \textbf{16}:6 (2017), 2357-2376.

\noindent [6] A. Kostianko and S. Zelik, ``Inertial manifolds for
1D reaction-diffusion-advection

systems. Part II: Periodic boundary conditions'', \textit{Comm.
Pure Appl. Anal}., \textbf{17}:1

(2018), 285-317.

\noindent [7] J. C. Robinson, \textit{Infinite-Dimensional
Dynamical Systems},  Cambridge Texts in Applied

Mathematics, Cambridge University Press, 2001.

\noindent [8] A. V. Romanov, ``Finite-dimensional limit dynamics
of dissipative parabolic equations'',

\textit{Sb. Mathematics}, \textbf{191}:3 (2000), 415-429.

\noindent [9] A. V. Romanov, ``Finite-dimensionality of dynamics
on an attractor for non-linear

parabolic equations'', \textit{Izvestia$:$ Mathematics},
\textbf{65}:5 (2001), 977-1001.

\noindent [10] A. V. Romanov, ``A Parabolic Equation with Nonlocal
Diffusion without a Smooth

Inertial Manifold'', \textit{Math. Notes}, 2014, \textbf{96}:4
(2014), 548-555.

\noindent [11] A. V. Romanov, ``Final dynamics of systems of
nonlinear parabolic equations

on the circle,'' \textit{AIMS Mathematics}, \textbf{6}:12 (2021),
13407-13422.

\noindent [12] R. Temam, \textit{Infinite-dimensional Dynamical
Systems in Mechanics and Physics},

Appl. Math. Sci., vol. \textbf{68} (2-nd ed.), Springer, N.Y.,
1997.

\noindent [13] H. Triebel, ``Theory of function spaces'', Monogr.
in Math., Vol. \textbf{78}, Birkhauser

Verlag, Basel--Boston--Stuttgart, 1983.

\noindent [14] S. Zelik, ``Inertial manifolds and
finite-dimensional reduction for dissipative PDEs,''

\textit{Proc. Roy. Soc. Edinburgh, Ser. A}, \textbf{144}:6 (2014),
1245-1327.

\bigskip \noindent School of Applied Mathematics,

\noindent National Research University Higher School of Economics,

\noindent 34 Tallinskaya St., Moscow, 123458 Russia

\noindent E-mail adress: av.romanov@hse.ru
\medskip

\end{document}